\documentclass{proc-l}
\usepackage{amsmath}
\usepackage{amsfonts}
\usepackage{amsthm}
\usepackage{amssymb}
\newtheorem{theorem}{Theorem}
\newtheorem*{theorem*}{Theorem}
\newtheorem*{acknowledgement*}{Acknowledgement}

\newtheorem{corollary}[theorem]{Corollary}

\newtheorem{lemma}[theorem]{Lemma}

\newtheorem{remark}[theorem]{Remark}

\newcommand{\Mn}[0]{\mathcal{M}^{n}}
\newcommand{\R}[1]{{\mathbb{R}}^{#1}}
\newcommand{\pd}[2]{\frac{\partial #1}{\partial#2}}

\newcommand{\pdt}[0]{\frac{\partial}{\partial t}}
\newcommand{\normsq}[1]{\left|#1\right|^2}
\newcommand{\heat}[0]{\left(\pdt -\Delta\right)}
\newcommand{\hess}[1]{\nabla\nabla #1}
\newcommand{\gsq}[1]{\normsq{\nabla #1}}
\newcommand{\hsq}[1]{\normsq{\hess{#1}}}
\newcommand{\norm}[1]{\left|#1\right|}
\newcommand{\ue}[0]{u_{\epsilon}}
\newcommand{\Me}[0]{M_{\epsilon}}

\newenvironment{mainproof}[1][Proof of Theorem 1]{\noindent\textbf{#1.} }{\ \rule{0.5em}{0.5em}}
\begin{document}

\title[Hamilton's estimate on complete manifolds]{Hamilton's gradient estimate for the heat kernel on complete manifolds}

\author[B. Kotschwar]{Brett L. Kotschwar}
\address{University of California, San Diego}
\email{bkotschw@math.ucsd.edu}
\begin{abstract}
	In this paper we extend a gradient estimate of R. Hamilton 
for positive solutions to the heat equation on closed manifolds
to bounded positive solutions on complete, non-compact manifolds 
with $Rc \geq -Kg$.  We accomplish this extension via
a maximum principle of L. Karp and P. Li and a Bernstein-type estimate
on the gradient of the solution.  An application of our result, 
together with the bounds of P. Li and S.T. Yau, yields an 
estimate on the gradient of the heat kernel for complete manifolds
with non-negative Ricci curvature that is sharp in the order
of $t$ for the heat kernel on ${\mathbb{R}}^n$.
\end{abstract}
\maketitle

\section{Introduction}
In \cite{Ham}, Richard Hamilton established the following estimate on
the gradient of the logarithm of a positive solution to the heat equation.

\begin{theorem*}{(Hamilton)}
Suppose $(\Mn, g)$ is a closed Riemannian manifold and $u$ a positive
solution to the heat equation on $\Mn$.  If $M > 0$ and $K \geq 0$ are
constants such that $Rc \geq -Kg$ and $u(x,t) \leq M$, then for all
$x\in\Mn$ and $t > 0$, one has
\begin{equation}
\label{eq:ham}
		t\normsq{\nabla \log(u)} \leq (1 + 2Kt)\log(M/u).
\end{equation}
\end{theorem*}

In this paper, we provide a proof that Hamilton's theorem also 
holds for complete, non-compact manifolds with Ricci curvature bounded
below.  Under the additional restriction of non-negative Ricci curvature, 
we then obtain, via the well-known bounds of Li and Yau \cite{LY},
the following estimate on the heat kernel.  (Recall that on a complete, non-compact manifold,
the heat kernel may be defined as the smallest positive fundamental solution
to the heat equation.)

\begin{theorem} 
\label{th:main}
Suppose $\Mn$ is a complete, non-compact manifold with $Rc \geq 0$,
and $H$ its heat kernel.  Then, for all $\delta > 0$,
there exists a constant $C = C(n, \delta)$ such that
\begin{equation}
\label{eq:kernelest}
	\normsq{\nabla \log(H(x,y,t))}\leq 
		\frac{2}{t}\left(C + \frac{d(x,y)^2}{(4-\delta)t}\right)
\end{equation}
for all $x$, $y$ in $\Mn$ and $t > 0$.
\end{theorem}

Theorem \ref{th:main} is sharp in the order of $t$ for the heat kernel
on $\R{n}$ and should be compared to the recent estimate of Souplet and Zhang
in \cite{SZ} which applies to the heat kernel
on manifolds with  $Rc(g)\geq -Kg$.  
In the special case $K=0$, inequality (\ref{eq:kernelest})
is comparable to their estimate at scales $d^2(x,y)\leq ct$ and offers an
improvement at scales $t << d^2(x,y)$.

Additionally, such an estimate is required to prove
that the integrand in the entropy functional 
$\mathcal{W}$ for the linear heat equation
(cf. \cite{Ni})
is pointwise non-positive for the fundamental solution to the heat equation.  
The proof that the integrand in 
Perelman's $\mathcal{W}$-functional is non-positive
for fundamental solutions to the conjugate heat equation seems also
to require a non-linear analog of this result (see \cite{N2}). Perhaps
the approach detailed here could serve as a model for the proof of such
an estimate. 

In Section \ref{sec:bern}, we obtain a Bernstein-type estimate  
for bounded solutions to the heat equation which affords pointwise control
on the product of the squared norm of the gradient by the elapsed time.  
The estimate is similar in form to those found by W.-X. Shi \cite{Shi} in
the Ricci Flow setting for derivatives of curvature.  
Such estimates have found considerable service in that setting
and continue to have importance in work towards the verification of the claims
of Perelman.

\begin{acknowledgement*}
This work is based on a suggestion of Professor Lei Ni, to whom 
the author
is also indebted for several useful discussions.  The author
would also like to thank Professor Bennett Chow for his support and encouragement,
and the referee for helpful suggestions pertaining to the exposition of this
paper.
\end{acknowledgement*}

\section{A Bernstein-Type Estimate}
\label{sec:bern}
Henceforth we shall assume that $(\Mn, g)$ is a smooth, complete, non-compact
Riemannian manifold with Ricci curvature uniformly bounded below by $-K$,
and for this section, suppose $u$ is a smooth solution to the heat equation on 
some open $U\subset\Mn$ for $0\leq t\leq T \leq \infty$, satisfying 
$|u|\leq M$.  Our aim is to establish a preliminary estimate on $\normsq{\nabla u}$ so that
the maximum principle of Ni and Tam \cite{NT} may be applied to the quantity
of interest in Hamilton's gradient estimate. 

To do this,  we employ a technique of 
W.-X. Shi \cite{Shi} from the estimation of derivatives of curvature under
the Ricci Flow (see also the treatment in the forthcoming book \cite{CNLRF}),
and define $F(x,t) = (4M^2 + u^2)\normsq{\nabla u}$ for $t > 0$.
The evolution of $F$ then possesses an advantageous
$-F^2$ term, as we see below. 
\begin{lemma}
\label{lm:F}
There exist positive constants $C_1$ and $C_2$ such that
\begin{equation}
	\pd{F}{t} \leq \Delta F +C_1KF - \frac{C_2}{M^4}F^2.
\end{equation}
\end{lemma}
\begin{proof}
We have
\[
	\heat \gsq{u} = -2\normsq{\nabla\nabla u} - 2Rc(\nabla u, \nabla u)
		\leq -2\normsq{\nabla\nabla u} + 2K\normsq{\nabla u}, 
\]
and
\[
	\heat u^2 = - 2\normsq{\nabla u}.
\]

So
\begin{align*}
\begin{split}
	\heat F &\leq (4M^2 +u^2)(2K\normsq{\nabla u}-2\normsq{\nabla\nabla u})
			-2\left|\nabla u \right|^4\\ 
		&\phantom{\leq (4M^2 +u^2)}	- 8u(\nabla\nabla u)(\nabla u, \nabla u)
\end{split}\\
\begin{split}
		&\leq -10u^2\hsq{u} + 10M^2K\gsq{u}
			- 2\norm{\nabla u}^4\\
		 &\phantom{-10M^2\hsq{u}}-8u(\hess u)(\nabla u, \nabla u). 
\end{split}
\end{align*}
Since
\[
		8|u|\norm{\hess{u}}\gsq{u} 
		\leq 10u^2\hsq{u} + \frac{8}{5}\norm{\nabla u}^4,
\]
and $4M^2 \gsq{u}\leq F \leq 5M^2\norm{\nabla u}^2$, we find 
\begin{align*}
	\heat F &\leq 10M^2K\gsq{u} 
		-\frac{2}{5}\norm{\nabla u}^4\\
		&\leq \frac{5}{2}KF -\frac{2}{625M^4}F^2
\end{align*}
as claimed.
\end{proof}

Now, as in \cite{LY}, for any $p\in \Mn$ and $R>0$ we may find a cut-off function
$\eta(x) = \eta_{p,R}(x)$ equal to $1$ on $B_p(R)$ and supported in $B_p(2R)$ 
satisfying the conditions
\begin{equation}
\label{eq:cutoffgrad}
	\gsq{\eta}\leq \frac{C_3}{R^2}\eta\\
\end{equation}
and
\begin{equation}
\label{eq:cutofflaplacian}
	\Delta \eta \geq -\frac{C_3}{R^2}\left(1 + R\sqrt{K}\right)
\end{equation}
for some $C_3 = C_3(n) > 0$.
Strictly speaking, the above estimates need only hold 
away from the cut-locus of $p$, however,
for the purposes of applying the maximum principle to $\eta F$, the well-known argument of Calabi \cite{Cal} allows us to assume that they hold everywhere.

The main result of this section is the following local estimate:

\begin{theorem}
	Suppose $u$ is a smooth solution to the heat equation satisfying
	$|u|\leq M$
	on $B_p(2R)\times [0,T]$ for some $p\in \Mn$ and $M$, $R$,  $T > 0$.
	Then there exists a constant
	$C_4 = C_4(K)$ such that 
\begin{equation}
\label{eq:bern}
	t\gsq{u} \leq C_4M^2\left(1 + T\left(1 + \frac{1}{R^2}\right)\right)
\end{equation}
holds on $B_p(R)\times [0,T]$
\end{theorem}
\begin{proof}
	Define $F$ as in Lemma \ref{lm:F}.   On $\mathrm{supp}(\eta)\times [0,T]$, we have, by Lemma \ref{lm:F} 
	and equations (\ref{eq:cutoffgrad}) and (\ref{eq:cutofflaplacian}),
\begin{align*}
	\heat(t\eta F) &= \eta F + t\eta\heat F - tF\Delta\eta 
		-2t\left\langle\nabla \eta, \nabla F\right\rangle\\
\begin{split}
	&= \left(\eta + 2t\frac{\gsq{\eta}}{\eta} - t\Delta \eta\right)F
	 	+ t\eta\heat F\\
	&\phantom{= (\eta} -2\left\langle\nabla(t\eta F), 
			\frac{\nabla\eta}{\eta}\right\rangle
\end{split}\\
\begin{split}
	&\leq \left(\left(1 +C_1Kt\right)\eta +
		 3t\frac{C_3}{R^2}\left(1+ R\sqrt{K}\right)\right)F
		 -\frac{C_2}{M^4}t\eta F^2\\
	&\phantom{\leq ((1}  -2\left\langle\nabla(t\eta F), 
			\frac{\nabla\eta}{\eta}\right\rangle.
\end{split}
\end{align*}	

If  $\eta F$ is not identically zero (i.e., if $u$ is not constant
on ${\mathrm{supp}}(\eta)$), then 
$t\eta F$ attains a positive maximum at $(x_0, t_0)\in \Mn\times (0,T]$.

At this point,
\[
	\nabla(t\eta F) = 0
\]
and
\[
 	\heat (t\eta F) \geq 0,
\]
so
\[
	\frac{C_2}{M^4}t_0\eta F^2 \leq  \left(1 +C_1KT +
		 3T\frac{C_3}{R^2}\left(1+ R\sqrt{K}\right)\right)F.
\]
Consequently, for any $(x,t)\in B_p(R)\times[0,T]$,
\begin{align*}
	tF(x,t) &= tF(x,t)\eta(x)\\ 
	&\leq t_0 F(x_0,t_0)\eta(x_0) \\
	&\leq \frac{M^4}{C_2}\left(1+ C_1KT + 3T\frac{C_3}{R^2}
			\left(1+R\sqrt{K}\right)\right).
\end{align*}
But $\gsq{u} \leq (1/4M^2)F$, and the claim follows.
\end{proof}

\begin{remark}  If $Rc\geq 0$, the above proof shows
\[
	t\gsq{u} \leq C M^2\left(1 + \frac{T}{R^2}\right)
\]
on $B_p(R)\times[0,T]$.
\end{remark}
Sending $R\to\infty$ in the penultimate line of the above proof, 
we at once obtain
\begin{corollary}
	Suppose the solution $u$ is defined on all of $\Mn\times[0,T]$. 
	Then there exists a constant $C_5$ such that
\begin{equation}
\label{eq:berncor}
	t\gsq{u} \leq C_5M^2(1+KT)
\end{equation}
on $\Mn\times[0,T]$.
\end{corollary}
\section{Proof of Main Theorem}

Next, using the estimate of the previous section, 
we apply a maximum principle due originally to 
Karp and Li, whose statement we found (in more generalized
form) in a paper of Ni and Tam.  
The statement of their theorem in our
(stationary metric) case is as follows.  Here 
$f_{+}(x,t) := \max\{f(x,t), 0\}$.
\begin{theorem*}{(Karp-Li, \cite{KL}; Ni-Tam, \cite{NT}, 1.2)}
	Suppose $(\Mn, g)$ is a complete Riemannian manifold and $f(x,t)$
	a smooth function on $\Mn\times[0,T]$ 
	such that $\heat f(x,t) \leq 0$ whenever
	$f(x,t)\leq 0$.  Assume that
\begin{equation}
\label{eq:mpcond}
	\int_0^T\int_{\Mn}\,e^{-ar^2(x)}f_+^2(x,s)\,dV ds < \infty
\end{equation}
for some $a > 0$, where $r(x)$ is the distance to $x$ from some fixed 
$p\in \Mn$.  If $f(x,0)\leq 0$ for all $x\in \Mn$, then $f(x,t)\leq 0$
for all $(x,t)\in \Mn\times[0,T]$.
\end{theorem*}

Hamilton's theorem in the complete case reads as 

\begin{theorem}
\label{th:ham}
  Suppose $(\Mn,g)$ is a complete manifold with
$Rc\geq -Kg$ for some $K\geq 0$. If $0 < u(x,t) \leq M$ is a solution
to the heat equation on $\Mn\times[0,T]$ for $0 < T \leq \infty$, then
\[
	t\gsq {\log u} \leq (1+ 2Kt)\log(M/u).
\]
\end{theorem}
\begin{proof}
	Defining $\ue = u + \epsilon$ for $\epsilon > 0$, we
	obtain a solution satisfying
	 $\epsilon < \ue \leq M + \epsilon :=\Me$.
	Once the estimate has been proved for $\ue$, the theorem will
	follow by letting $\epsilon\to 0$.

As in \cite{Ham}, the function
\[
	P(x,t) := \varphi(t)\frac{\gsq{\ue}}{\ue}-\ue\log(\frac{\Me}{\ue}),
\]
where $\varphi(t) := t/(1+2Kt)$, satisfies $\heat P(x,t)\leq 0$ and
\[
	P(x,0) = -\ue\log(\Me/\ue) \leq 0.
\]

By our assumptions on $\ue$, we also have
\[
	P_{+}(x,t) \leq \frac{1}{\epsilon}t\gsq{\ue}.
\]
Thus using equation (\ref{eq:berncor}), for any $p\in \Mn$, and $R > 0$, we have
\begin{equation*}
\begin{split}
	\int_0^T\int_{B_p(R)}\,e^{-r^2(x)}P_{+}^2(x,t)\,dV dt 
		\leq\frac{1}{\epsilon^2}\int_0^T\int_{B_p(R)}\,e^{-r^2(x)}
			(t\gsq{\ue})^2\,dVdt\\
		\phantom{\int_0^T\int_{B_p(R)}\,}\leq \frac{C_5^2\Me^4}{\epsilon^2}
		\left(1 + KT\right)^2 
		\int_0^T\int_{\Mn}\,e^{-r^2(x)}\,dV dt.
\end{split}
\end{equation*}
Since we assume that $Rc \geq -Kg$, it follows from the Bishop volume comparison
theorem that the rightmost integral in the above inequality is finite.

Hence,
\begin{align*}
	\int_0^T\int_{\Mn}\,e^{-r^2(x)}P_{+}^2(x,t)\,dV dt &\leq
	\liminf_{R\to\infty}\int_0^T\int_{B_p(R)}\,e^{-r^2(x)}P_{+}^2(x,t)\,dV dt\\
	&<\infty,
\end{align*}
and we conclude that $P(x,t) \leq 0$ for all $t\leq T$.
\end{proof}

\begin{mainproof}
Let $H(x,y,t)$ denote the heat kernel of $(\Mn, g)$.  For any $t>0$ 
and $y\in \Mn$, set
$u(x,s) := H(x,y,s +t/2)$.  Then $u$ is a smooth, positive solution to the heat equation 
on 
$[0,\infty)$.  By Corollary 3.1 and Theorem 4.1 of \cite{LY}, for all 
$\delta > 0$, there is a constant $C_6 = C_6(\delta) > 0$ such
that $u$ satisfies
\begin{equation}
\label{eq:ly}
	C_6^{-1}V\left(\sqrt{s+t/2}\right)^{-1}e^{\frac{-d^2(x,y)}
		{(4-\delta)(s+t/2)}}
	\leq u(x,s)\leq C_6V\left(\sqrt{s+t/2}\right)^{-1} 
\end{equation}
for all $x$, $y\in \Mn$, and $s \geq 0$, where 
$V(\sqrt{s+t/2}):= {\mathrm{Vol}}(B_y(\sqrt{t +s/2}))$.  

Defining $M = C_6 V\left(\sqrt{t/2}\right)^{-1}$,
the left inequality in (\ref{eq:ly}) implies $u \leq M$ for all $x$ and $s$.
Moreover, since we assume $Rc\geq 0$, there exists a positive constant 
$C_7 = C_7(n)$ such that for all $0\leq s\leq t/2$
\[
	V\left(\sqrt{t/2 + s}\right) \leq V\left(\sqrt{t}\right) 
	\leq C_7 V\left(\sqrt{t/2}\right). 
\]
Thus, by the right-hand inequality in (\ref{eq:ly}) and Theorem \ref{th:ham}, 
we have
\[
	s\gsq{\log u}\leq \log\left(\frac{M}{u}\right) 
 \leq 
\left(\log (C_6^2C_7) +\frac{d^2(x,y)}
{(4-\delta)(s+t/2)}\right)
\]
on $\Mn\times [0, t/2]$.  \

Setting 
$C=\log (C_6^2(\delta)C_7(n))$ and evaluating at $s = t/2$, we conclude that
\[
	(t/2)\gsq{\log H}(x,y,t) = (t/2)\gsq{u}(x,t/2) \leq 
\left(C + \frac{d^2(x,y)}
{(4-\delta)t}\right)
\]
 for all $x$, $y \in \Mn$ and $t >0$.
\end{mainproof}

\end{document}